\documentclass[11pt,fleqn]{article}

\usepackage{amssymb,latexsym,amsmath, amsthm}
\usepackage{times}

\newcommand{\pr}[1]{{\rm {\mathsf P}}\{#1\}}

\newcommand{\expect}[1]{{\rm {\mathsf E}}\{#1\}}

\theoremstyle{plain}
\newtheorem{thm}{Theorem}[section]
\newtheorem{cor}{Corollary}[section]
\newtheorem{lem}{Lemma}[section]

\theoremstyle{definition}
\newtheorem{defn}{Definition}[section]

\theoremstyle{remark}
\newtheorem{rem}{Remark}

\begin{document}
\section*{\centering\Large\rm
          Asymptotics for first-passage times of L\'evy processes and random walks
}
\section*{\centering\large\rm
D. Denisov,\footnote{Address: School of MACS,
Heriot-Watt University,
Edinburgh EH14 4AS, UK.
E-mail address: Denisov@ma.hw.ac.uk}
and V. Shneer,\footnote{Address: EURANDOM,
Eindhoven University of Technology,  P.O. Box 513 - 5600 MB Eindhoven, The Netherlands.
E-mail address: Shneer@eurandom.tue.nl}}
\section*{\centering\normalsize\it
Eurandom and  Heriot-Watt University}

\begin{abstract}
We study the exact asymptotics for the distribution of the first
time $\tau_x$ a L\'evy process $X_t$ crosses  a negative level $-x$.
We prove  that $\mathbf P(\tau_x>t)\sim V(x)\mathbf P(X_t\ge 0)/t$
as $t\to\infty$ for a certain function $V(x)$. Using known results
for the large deviations of random walks we obtain asymptotics for
$\mathbf P(\tau_x>t)$ explicitly in both light and heavy tailed
cases. We also apply our results to find asymptotics for the
distribution of the busy period in an M/G/1 queue.\\
{\it Keywords:}
L\'evy processes, random walks, busy period, first-passage times, subexponential distributions,
 large deviations,single server queue.\\
{\it AMS Subject Classification.} Primary: 60G50, 60G51; secondary 60K25.

\end{abstract}


\section*{Introduction}

Let $\{X_t\}_{t\ge 0}$ be a L\'evy process with the characteristic
function $\expect{e^{i\theta X_t}}=e^{t\Psi(\theta)}, t\ge 0$, where
$\Psi$ is given by the L\'evy-Khinchine formula\cite{Sato}
\begin{equation}\label{levi-khinchine}
\Psi(\theta)=iA\theta-\frac{1}{2}\sigma^2\theta^2+
\int_{-\infty}^{+\infty}(e^{i\theta x}-1-i\theta x {\bf
1}_{[-1,1]})\Pi ({\rm dx}). \end{equation} For $x\ge 0$ let
$$
\tau_x=\min\{t\ge 0: X_t<-x\}
$$
be the first-passage time. Throughout we assume that the L\'evy
process $X_t$ drifts to $-\infty$ a.s. Rogozin's criterion
\cite{rogozin} (see also \cite[p. 167]{bertoin} or \cite[theorem
48.1]{Sato}) says that $X_t\to-\infty$ if and only if
\begin{eqnarray}\label{trans}
\int_1^{\infty}t^{-1}\pr{X_t\ge 0}{\rm dt}<\infty.
\end{eqnarray}
This assumption implies that $\tau_x$ is a proper random variable
with finite expectation:
$$
\tau_x<\infty,\quad \expect{\tau_x}<\infty.
$$
The aim of this work is to study asymptotics
$$
\pr{\tau_x>t}
$$
when $x > 0$ is fixed and $t \to \infty$.

An analogue of the L\'evy processes for discrete time is random
walks.  Let $S_n=\xi_1+\cdots+\xi_n$ be a random walk with i.i.d.
increments and assume that $S_n \to -\infty$ a.s.
Then an analogue of $\tau_x$ is the stopping time
$\nu_x=\min\{n\ge 1:S_n< -x\}$. Asymptotics for $\nu_x$ have been
studied in \cite{eh} for $x=0$, and in \cite{doney},
\cite{bertoin_doney} for $x>0$. In \cite{eh}, it is shown that
$\pr{\nu_0>n}\sim\pr{S_n\ge 0}/n$ if the latter is a
subexponential sequence (see definition below).   In \cite{doney}
and \cite{bertoin_doney}, the asymptotics for $\nu_x$ have been
found when $x>0$. In these papers the authors considered
separately several classes of distribution of $\xi$: regularly
varying distributions, distributions satisfying Cram\'er's
condition and light-tailed distributions which do not satisfy
Cram\'er's condition. For each of these classes, they show that
$\pr{\nu_x>n}$ is asymptotically proportional to $\pr{S_n\ge 0}/n,$
using large-deviations asymptotics for $\pr{S_n\ge 0}$. In
\cite{borovkov2} (see also \cite{borovkov}) it has been shown that
the same asymptotic equivalence holds for the so-called
semi-exponential distributions with parameter smaller that $1/2$
(see definition in \cite{borovkov}). In \cite{Mogulskii} the
latter result has been generalised to the case of semi-exponential
distributions with parameter smaller than $1$. We note also that
in \cite{Mogulskii} $x$ may depend on $n$. We also mention the
paper \cite{borovkov3} where the same problem was studied under
Cram\'er's assumptions and in the case when $x$ may depend on $n$.

For L\'evy processes, the asymptotics for $\tau_x$ have been
obtained in \cite{KR} for the light-tailed distributions.
Necessary and sufficient conditions for the existence of moments
of $\tau_x$ have been obtained in \cite{doney_maller}.

In our work we develop the approach proposed in \cite{doney},
\cite{bertoin_doney}. The main results of our paper are Theorems
\ref{thm1} and \ref{thm_disc}. Theorem \ref{thm1} states that
under some natural assumptions,
\begin{eqnarray*}
\pr{\tau_x>t}&\sim & V(x)\frac{\pr{X_t\ge 0}}{t},
\end{eqnarray*}
for some  function $V(x)$ depending only on $x$. Theorem
\ref{thm_disc} shows that under identical conditions for both L\'evy
processes and random walks,
\begin{eqnarray*}
\pr{\nu_x>n}&\sim &V _{rw}(x)\frac{\pr{S_n\ge 0}}{n},\\
\pr{\tau_x>t}&\sim & V(x) e^{-\gamma(t-[t])}\frac{\pr{X_{[t]}\ge
0}}{t}
\end{eqnarray*}
for some $\gamma \ge 0$, where by $[t]$ we denote the integer part
of $t$ (the largest integer smaller than $t$) and functions $V(x)$
and $V_{rw}(x)$ depend only on $x$. It is worth mentioning that for
the case of L\'evy processes both conditions of Theorem
\ref{thm_disc} and its result are given in terms of the values of
the process at discrete times. Therefore, the problem of finding
asymptotics for $\pr{\tau_x > t}$ and $\pr{\nu_x>n}$ is reduced to
finding asymptotics for $\pr{S_n \ge 0}$, or $\pr{\tilde{S}_n \ge
na}$ where $a = - \expect{\xi_1}$ and $\tilde{S}_n = S_n + na$ is a
random walk with zero drift. This is a problem of large deviations
of sums of independent identically distributed random variables
which is extensively studied in literature. We apply known results
to obtain explicit asymptotics in various cases. It appears that in
all cases when asymptotics for $\pr{S_n > 0}$ can be found
explicitly, the conditions of Theorem \ref{thm_disc} are satisfied
and hence, asymptotics for the tail distribution of $\tau_x$ and
$\nu_x$ can also be found explicitly.

We consider distributions with heavy tails (such that
$\expect{e^{\varepsilon \xi_1}} = \infty$ for all $\varepsilon >
0$) and distributions with light tails (for which the latter
condition is not fulfilled) separately. Moreover, it has been pointed out by
various authors that for the problem of large deviations of sums
of heavy-tailed random variables, one should indicate two classes
of distributions: those with tails lighter and heavier than
$e^{-\sqrt t}$ (we say that the tail of a distribution $F$ is
lighter than a function $f$ if $\overline F(t)/f(t) \to 0$ as $t
\to \infty$ and heavier than a function $f$ if $\overline F(t)/f(t)
\to \infty$). For the case of heavy-tailed distributions some new results on large deviations are presented in
\cite{DDS}. We also state some of these results in the present paper,
however, we do not concentrate on the problem of large deviations
here.

Theorem \ref{thm_disc} can also be applied to the case of
light-tailed distributions under some further assumptions. In
particular, the conditions of Theorem \ref{thm_disc} are fulfilled
if the distribution of $\xi_1$ satisfies the so-called Cram\'er's
(or classical) conditions. With the help of Theorem \ref{thm_disc}
we can also cover the so-called intermediate case, when $\xi_1$ has
a distribution with a light tail but does not satisfy Cram\'er's
condition.

Another motivation for our work was to find asymptotics for the
busy period in a stable $M/G/1$ queue. Let  $A_1, A_2, \ldots$
and $B_1, B_2, \ldots$ be two mutually independent sequences each
consisting of independent and identically distributed random
variables. Assume that $\{A_i\}$ are inter-arrival times and
$\{B_i\}$ are service times. We assume throughout that
$\expect{B_1}/\expect{A_1} = \rho < 1$ so that the system is
stable. We use the common notation: we denote by $M/G/1$ the system
when $A_i$ are exponential random variables; in the case of a
general
i.i.d. sequence $\{A_i\}$, we denote the system by $GI/GI/1$. 
Denote $N(t) = \max\{n: A_1 + \ldots + A_n \le t\}$. Put $X_0 = 0$
and
\begin{equation} \label{eq_bp}
X_t = \sum\limits_{i=1}^{N(t)} B_i-t.
\end{equation}
Then the busy period of the system with initial work $x>0$ may be
defined as
\[
bp(x) = \inf\{t: X_t < -x\}.
\]
Hence, in an $M/G/1$ queue, finding asymptotics for the tail of
$bp(x)$ is equivalent to finding asymptotics for the tail of
$\tau_x$ when $X_t$ is a compound Poisson process without negative
jumps.

The tail behavior of the busy period in these systems has been
studied by various authors under different assumptions. Under
Cram\'er-type assumptions, asymptotics for the $M/G/1$ setting  were
studied in \cite{AW} and for the $GI/G/1$ setting --- in \cite{PR}.
Most of the papers on the tail behaviour of the busy period are
devoted to studying the case when $B_1$ has a subexponential
distribution. All these papers investigate the asymptotic behaviour
of $bp$ --- busy period of the queue under the condition that the
first customer arriving at the system finds it empty. In \cite{MT},
it was shown that if $B_1$ has a regularly varying distribution then
\begin{equation} \label{asymptotics}
\pr{bp > t} \sim \expect{\nu_0} \pr{B_1 > (1-\rho) t}
\end{equation}
as $t \to \infty$. This result has been generalised in
\cite{Zwart} to the case of a $GI/G/1$ queue and under the
assumption that the tail $\overline B(t) = \pr{B_1 > t}$ satisfies
an extended regular variation condition (see \cite{Bingham}).

Later on, it has been shown in \cite{BDK} and \cite{JM}, that the
asymptotics (\ref{asymptotics}) hold for the $GI/G/1$ model for
another subclass of heavy-tailed distributions which includes the
Weibull distributions with parameter $\alpha<1/2.$  The tails of
the distributions considered in \cite{BDK} and \cite{JM} are
heavier than $e^{-\sqrt t}$. As is shown in \cite{AKS} (see also
\cite{foss_korshunov}), the latter condition is  crucial for the
asymptotics (\ref{asymptotics}) to hold.

The method proposed in this paper allows to find asymptotics for
$\pr{bp(x)> t}$ in the $M/G/1$ queue for both  light- and
heavy-tailed distributions. Moreover, we are able to obtain these
asymptotics when $\pr{B>t}$ is lighter than $e^{-\sqrt t}$ but
still heavier than any exponential distribution. Using the results
on tail asymptotics of the distribution of $bp(x)$, we can also
obtain the results for the tail asymptotics of the distribution of
$bp$.


The paper is organised as follows. In section \ref{sec_1} we present
Theorems \ref{thm1} and \ref{thm_disc} that reduce the problem of
finding asymptotics for $\pr{\tau_x>t}$ and $\pr{\nu_x>n}$ to
studying asymptotics of $\pr{S_n \ge  0}$. In section \ref{sec_2} we
consider 4 classes of distributions: heavy-tailed distributions I
(with tails heavier than $e^{-\sqrt t}$), heavy-tailed distributions
II (with tails lighter than $e^{-\sqrt t}$), distributions
satisfying Cram\'er's condition and distributions forming an
intermediate case (distributions with light tails not satisfying
Cram\'er's condition). For each of these cases we give known results
on asymptotics of $\pr{S_n \ge 0}$, show that the conditions of
Theorem \ref{thm_disc} are satisfied and hence, obtain results on
the tail asymptotics of the distributions of $\tau_x$ and $\nu_x$.
Appendix A is devoted to the proofs of Theorems \ref{thm1} and
\ref{thm_disc}, in Appendix B we present some known results on
L\'evy processes that are used in our paper.

\section{Main results} \label{sec_1}

In this section we present Theorems \ref{thm1} and \ref{thm_disc}
which connect the asymptotics for $\tau_x$ and $\nu_x$ with the
asymptotics for $\pr{S_n \ge 0}$. Before stating general Theorems,
we need to introduce
\begin{defn}
A function $f:\mathbf R^+\to \mathbf R^+$ belongs to the class
$\mathcal Sd(\gamma)$ with $\gamma\ge 0$ if, starting from some
moment $t$, $f(t)>0$ and
\begin{eqnarray}\label{lt}
\lim_{t\to\infty}f(t-y)/f(t)=e^{\gamma y},\quad  y\in \mathbf
R;\end{eqnarray}
\begin{eqnarray}\label{subex}
\lim_{t\to\infty}\frac{f^{*2}(t)}{f(t)}=\lim_{t\to\infty}
\frac{\int_{0}^{t} f(t-y)f(y){\rm dy}}{f(t)}= 2d =2\int_0^\infty
e^{\gamma y}f(y){\rm dy}.
\end{eqnarray}
The class $\mathcal Sd:=\mathcal Sd(0)$ is called the class of
subexponential densities.
\end{defn}
A discrete-time analogue of this definition is
\begin{defn}
A sequence $\{a_n\}_{n\ge 0}$ belongs to the class $\mathcal
Ss(\gamma)$ with $\gamma\ge 0$ if starting from some index $n, a_n>0$
and
\begin{eqnarray}\label{lt_disc}
\lim_{n\to\infty}a_{n-1}/a_n=e^{\gamma},\quad  \end{eqnarray}
\begin{eqnarray}\label{subex_disc}
\lim_{n\to\infty}\frac{a^{*2}_n}{a_n}\equiv\lim_{n\to\infty}\frac{\sum_{i=0}^n
a_i a_{n-i} }{a_n}= 2d =2\sum_{i=0}^\infty a_i e^{\gamma i}.
\end{eqnarray}
The class $\mathcal Ss:=\mathcal Ss(0)$ is called the class of
subexponential sequences.
\end{defn}

\begin{thm}\label{thm1}
Let the function
$$\frac{\pr{X_t\ge0}}{t},\quad t\ge 1$$ belong to the
class $\mathcal Sd(\gamma)$. In addition assume that for some
$\alpha\ge 0$
\begin{eqnarray}\label{cond1}
\lim_{t\to\infty}\frac{\pr{X_t\ge 0}}{\pr{X_t\ge y}}=e^{\alpha
y},\quad \mbox{ for any fixed } y.
\end{eqnarray}
Then,
\begin{eqnarray}
\pr{\tau_x>t}\sim V(x)\frac{\pr{X_t\ge 0}}{t},
\end{eqnarray}
for any $x$, a point of continuity of the function
\begin{eqnarray}\label{eq_cx}
V(x)\equiv\biggl\{
\begin{array}{cc}
  \expect{\tau_x},& \gamma=\alpha=0   \\
  e^{\alpha x}\int_0^\infty e^{\gamma t}\expect{e^{\alpha N_t};|N_t|\le x}{\rm dt},&\mbox{otherwise}
\end{array}
\end{eqnarray}
where $N_t = \inf\limits_{0 \le s \le t} X_s$.
\end{thm}

In the next Theorem we show that it is possible to obtain
asymptotics for $\tau_x$ and $\nu_x$ under the same conditions for
both random walks and L\'evy processes.
\begin{thm}\label{thm_disc}
Let $X_t$ be either a L\'evy process or a random walk. Let the
sequence
$$\frac{\pr{X_n\ge 0}}{n},\quad n\in \mathbf N$$ belong to the class
$\mathcal Ss(\gamma)$. In addition assume that for some $\alpha\ge
0$
\begin{eqnarray}\label{cond1_disc}
\lim_{n\to\infty}\frac{\pr{X_n\ge 0}}{\pr{X_n\ge y}}=e^{\alpha
y},\quad \mbox{ for any fixed } y,\quad n\in\mathbf N.
\end{eqnarray}
Then, if $X_t$ is a L\'evy process, for any $x$,
\begin{eqnarray*}
\pr{\tau_x>t}&\sim &V(x)e^{-\gamma(t-[t])}\frac{\pr{X_{[t]}\ge
0}}{t},
\end{eqnarray*}
where $V(x)$ is defined in Theorem \ref{thm1}.

If $X_n$ is a random walk, then
\begin{eqnarray*}
\pr{\tau_x>n}&\sim & V_{rw}(x)\frac{\pr{X_{n}\ge 0}}{n},
\end{eqnarray*}
where
\[
V_{rw}(x)=\biggl\{
\begin{array}{cc}
  \mathbf E\nu_x,& \gamma=\alpha=0   \\
  e^{\alpha x}\sum_{k=0}^\infty e^{\gamma k}\mathbf E\{e^{\alpha N_k};|N_k|\le x\},&\mbox{otherwise}
\end{array}
\]
where $N_k = \min\limits_{0 \le l \le k} X_l$.
\end{thm}

\begin{rem} \label{equality}
The conditions of Theorem \ref{thm1} and (or) conditions of Theorem
\ref{thm_disc} imply that $e^{-\gamma}=\expect{e^{\alpha X_1}}$. The
proof of this fact is given in Appendix A. Note also that this fact
implies that $\alpha=0$ if and only if $\gamma=0$. This corresponds
to the subexponential case.
\end{rem}

\section{Explicit results} \label{sec_2}

This section consists of 4 subsections. Each of these subsections is
devoted to a class of distributions for which we present known
results on large deviations of sums of random variables, with the
help of these results we show that the conditions of Theorem
\ref{thm_disc} are satisfied and as a result we obtain the
asymptotics for the tail distributions of $\tau_x$ and $\nu_x$.
First, we prove Theorem \ref{levy_less_half} in which we study the
case when $-\ln\pr{X_1>t}=o(\sqrt t).$ Further, in Theorem
\ref{st_greater_half}, we analyse the case when $-\ln\pr{X_1>t}$ is
regularly varying with parameter $\alpha\in[1/2,1)$. For the first
two cases, we use some results from paper \cite{DDS}. Then, in
Theorem \ref{cramers_case} we give the asymptotics for Cram\'er's
case. This includes (partially) the distributions with exponential
tails and tails which are lighter than exponential. Finally, in
Theorem \ref{thm_intermediate}, we analyse distributions with
exponential tails that are not covered by Cram\'er's case. As
corollaries we give corresponding results for the tail asymptotics
of the busy period of an $M/G/1$ queue.

\subsection{Heavy-tailed distributions I} \label{sec_21}





The following Theorem easily follows from \cite[Corollary 2.1]{DDS}.

\begin{thm} \label{ld_less_half}
Let $S_n = \sum\limits_{i=1}^n \xi_i$ be a random walk. Let
$\expect{\xi_1}=0$ and $\expect{|\xi_1|^{\kappa}} < \infty$ for
some $\kappa \in (1, 2]$. Assume that
\begin{equation} \label{insens_1}
 \frac{\overline F(n -
n^{1/\kappa})}{\overline F(n)} \to 1
\end{equation}
as $n \to \infty$ and assume also that $x\to x^\kappa\mathbf P(\xi_1>x)$ either belongs to
$\mathcal Sd(0)$ or is O-regularly varying.
 Then for any $a > 0$,
\begin{equation} \label{normal_asymptotics}
\pr{S_n > na} \sim n \pr{\xi_1 > na}
\end{equation}
as $n \to \infty$.
\end{thm}

For the tail asymptotics of $\tau_x$ and $\nu_x$ the following is
true.

\begin{thm} \label{levy_less_half}
Let $X_t$ be either a L\'evy process or a random walk.
Assume that the distribution of $X_1$ satisfies the conditions of
Theorem \ref{ld_less_half}. Let $\expect{X_1}=-a<0$. Then
\begin{eqnarray} \label{as_levy_less_half}
\pr{\tau_x > t} &\sim &\expect{\tau_x} \pr{X_1 > ta}\sim
\expect{\tau_x} \overline \Pi(ta),\quad t\to\infty; \\
\label{as_rw_less_half}\pr{\nu_x > n} &\sim &\expect{\nu_x}
\pr{X_1 > na},\quad n\to\infty.
\end{eqnarray}
\end{thm}

\begin{rem} Note that the conditions of the latter Theorem imply
that $\overline F(y-\sqrt y)\sim\overline F(y)$. It, in turn,
implies that $-\ln \overline F(y)=o(\sqrt y). $ Thus, we again
have the Weibull distribution with parameter $1/2$ as a boundary.
\end{rem}

\begin{rem}
The proof of Theorem \ref{levy_less_half} relies on Theorem
\ref{ld_less_half} which is a result of \cite{DDS}. However, the
same asymptotics for $\tau_x$ and $\nu_x$ may be obtained for any
distributions satisfying the asymptotic equivalence $\pr{S_n
> na} \sim n \pr{\xi_1 > na}$. For instance, from the results of
\cite{nagaev1969} it follows that such asymptotics hold for
regularly varying distributions. The results of \cite{rozovskii}
imply that the same holds for Weibull-type distributions with a
parameter smaller than $1/2$. In \cite{DDS} it is shown that
Theorem \ref{ld_less_half} includes all the results known
beforehand.
\end{rem}


{\sc Proof of Theorem \ref{levy_less_half}.} We should check the
conditions of Theorem \ref{thm_disc}. First, it follows from
Theorem \ref{ld_less_half} that $\pr{X_n\ge 0}\sim \pr{X_{n+1}\ge
0}$ and $\pr{X_n\ge y}\sim \pr{X_n\ge 0}$. It is then easy to
check that $\pr{X_n\ge 0}/n$ is a subexponential sequence.

Then, it follows from Theorem \ref{ld_less_half} and Theorem
\ref{thm_disc} that
$$
\pr{X_t\ge 0}\sim\pr{X_{[t]} \ge 0}= \pr{X_{[t]}+[t]a\ge [t]a}\sim
[t]\pr{X_1\ge [t]a}.
$$
Then,
$$
\pr{\tau_x > t} \sim \expect{\tau_x} \frac{\pr{X_t \ge 0}}{t}\sim
\expect{\tau_x} \pr{X_1\ge [t]a}\sim \expect{\tau_x} \pr{X_1>ta}.
$$
The proof of Theorem \ref{levy_less_half} is complete.

We also present a direct corollary of Theorem
\ref{levy_less_half}.

\begin{cor}\label{cor2}
Let $\expect{B_1^\kappa} < \infty$ and $\expect{A_1^\kappa} <
\infty$ for some $\kappa \in (1, 2]$. Assume that the distribution
of $B_1$ satisfies the conditions of  Theorem~\ref{ld_less_half}.
Then
\begin{equation} \label{as_bp_less_half}
\pr{bp(x) > t} \sim \frac{x}{\mathbf E A_1} \frac{1}{1-\rho} \pr{B_1
> (1-\rho) t}
\end{equation}
as $t \to \infty$ for any fixed $x>0$.
\end{cor}

{\sc Proof of Corollary \ref{cor2}.} In the case of an $M/G/1$
system, case $X_1 = \sum\limits_{i=1}^{N(1)} X_i - 1$ (see
(\ref{eq_bp})). Hence, $\expect{X_1}=\rho - 1$ and $\mathsf P(X_1
> t) \sim \dfrac{1}{\mathbf E A_1} \mathsf P(B_1 > t)$ as $t \to \infty$.
 Theorem
\ref{levy_less_half} yields that
\[
\pr{bp(x) > t} \sim \expect{bp(x)}\pr{B_1 > (1-\rho) t}
\]
as $t \to \infty$. It remains to note that $\mathbf E bp(x) =
\dfrac{x}{1-\rho}$ (this, for example, can be obtained from
\cite[p.261, (4.94)]{Cohen}). The proof of Corollary \ref{cor2} is complete.

\subsection{Heavy-tailed distributions II} \label{sec_22}

Denote $g(x) = -\ln \overline F(x)$. In this subsection we
consider the case
$$
\limsup\frac{g(x)}{\sqrt x}>0.
$$ For this we use
\cite[theorem 5a]{rozovskii}.
\begin{thm}(\cite[theorem 5a]{rozovskii}) \label{rw_greater_half}
Let $\pr{\xi_1 > y} \sim e^{-g(y)}$ as $y \to \infty$ with a
doubly differentiable function $g$ such that $g''(y)$ does not
decrease for $y \ge y_o$ and $y g''(y) \sim (\beta-1) g'(y)$ as $y
\to \infty$ for some $\beta \in (0,1)$. Let
\begin{equation} \label{def_k}
k = \max\left\{l \in \{1,2,\ldots\}: \limsup\limits_{z \to \infty}
\frac{g(z)}{z^{l/(l+1)}} > 0\right\}.
\end{equation}
Let $$\expect{\xi}=0,\expect{\xi^2}=1,
\expect{|\xi|^{k+3}}<\infty,
$$
\begin{equation} \label{def_Q_k}
R(y)=g(y)+\frac{(t-y)^2}{2n}-\sum_{i=1}^k\lambda_{i-1}\frac{(t-y)^{i+2}}{n^{i+1}}.
\end{equation}
Let  $y_*$  be the maximal solution of $R'(y) = 0$. Then $y^*\le
t-\sqrt n$ and
\begin{equation} \label{as_levy_greater_half}
\pr{S_n > t} \sim n\sqrt\frac{1}{nR''(y_*)}\exp{\{-R(y_*)\}},\quad
n\to\infty
\end{equation}
uniformly in $t>1.6 \eta(n)$, where $\eta(z)$ is such that
$\eta^2(z)/g(\eta(z))\sim z, z\to \infty$. Here, $\lambda_i$ are the
coefficients of the Cram\'er series (see \cite{petrov2} for the
definition).
\end{thm}

\begin{rem} Note that the conditions of Theorem~\ref{rw_greater_half} imply that $g''(y)$ is a
regularly varying function with parameter $(\beta-2)$. This fact
follows from the monotonicity of $g''$ and Karamata's Theorem.
Then,  $g'(y)$ is regularly varying with parameter $(\beta-1)$ and
$g(y)$ is regularly varying with parameter $\beta$. Also, under
these conditions $\eta(z)$ may be equivalently defined as a
function such that $|g''(\eta(z))|\sim \beta(1-\beta)/z, z\to
\infty$. Therefore, $\eta(z)$ is a monotone regularly varying
function with parameter $1/(2-\beta)$.
\end{rem}

\begin{rem} In the statement of \cite[theorem 5a]{rozovskii} it is not
said that $y^*\le t-\sqrt n$, but  one can find this assertion in
the proof of \cite[Lemma 3a]{rozovskii}.
\end{rem}

We find it difficult to apply Theorem \ref{rw_greater_half}
directly, since it gives the asymptotics in terms of the maximal
solution to an equation. Therefore, we use the approach developed in
\cite{foss_korshunov} to simplify this equation.
\begin{lem}\label{lem3}
Let all conditions of Theorem \ref{rw_greater_half} hold. Let
$t_n\to\infty$ be a sequence  such that $t_n\ge 1.6 \eta(n)$.  Let
 $y_n$ be any sequence such that $y_n\sim t_n$ and
\begin{eqnarray}\label{cond_1half}
R'(y_n)=o(1/\sqrt n).
\end{eqnarray}
Then,
$$
\pr{S_n > t_n} \sim
n\sqrt\frac{1}{nR''(y_n)}\exp{\{-R(y_n)\}},\quad n\to\infty.
$$
Also, for any sequence $\varepsilon_n=o(\sqrt n)$,
$$
\pr{S_n>t_n}\sim\pr{S_n>t_n+\varepsilon_n}.
$$
\end{lem}

{\sc Proof of lemma \ref{lem3}.}  First, we note that since $g''$
is monotone and regularly varying it is true that
$$
R'(y+z)-R'(y)=zR''(y)(1+o(1)), \quad y\to\infty, z=o(y).
$$
Also, $R''(y_n)=(g''(y_n)+1/n)(1+o(1))$ and
\begin{multline*}
|g''(y_n)|\le |g''(1+o(1)t_n)|\le(1+o(1))|g''(1.6\eta(n))|\\ \le
(1+o(1))|g''(\eta(n))|=(1+o(1))\beta(1-\beta)/n\le 1/(4n)
\end{multline*}

Then, for any $\varepsilon>0,$
\begin{eqnarray*}R'(y_n+\varepsilon\sqrt n)=R'(y_n)+\varepsilon\sqrt n
R''(y_n)\ge o(1/\sqrt n)+3/4\varepsilon/\sqrt n>0\\
R'(y_n-\varepsilon\sqrt n)=R'(y_n)-\varepsilon\sqrt n R''(y_n)\le
o(1/\sqrt n)-1/4\varepsilon/\sqrt n<0.
\end{eqnarray*}
Since $R'$ is continuous there exists a sequence
$\beta_n\in(y_n-o(\sqrt n),y_n+o(\sqrt n))$ such that
 $R'(\beta_n)=0$ and $\beta_n\sim t_n$.  Further, if there exists
some other solution $\beta'_n>\beta_n$, then with necessity
$\beta'_n\sim t_n\sim\beta_n$. But this is not possible since
$R''(y)$ is positive on the interval $(\beta_n, t_n)$.

To prove the first statement of the lemma, note
\begin{multline*}
R(y_n)-R(\beta_n)=R'(\beta_n)(\beta_n-y_n)+(1+o(1))R''(\beta_n)\frac{(\beta_n-y_n)^2}{2}\\
=(1+o(1))R''(\beta_n)\frac{(\beta_n-y_n)^2}{2}\sim O(\frac{1}{n})
o(\sqrt n)^2=o(1).
\end{multline*}
To prove the second statement of the lemma we should note that if
$\beta_n$ is a solution sequence for the equation $R_1(\beta_n)=0$
for the first sequence $t_n$, then the corresponding equation for
the sequence $t_n+\varepsilon_n$ is $R_2(\beta_n)=o(1/\sqrt n)$.
Then we should just apply the first statement of the lemma.
The proof of lemma \ref{lem3} is complete.

We shall now concentrate on the case $t_n = na$, the case needed
for our purposes.

\begin{cor}\label{cor1}
Let all conditions of Theorem \ref{rw_greater_half} hold. Let
$t_n=na$, where $a>0.$  Let
 $y_n$ be any sequence such that $y_n\sim t_n$ and
Condition (\ref{cond_1half}) holds. Then,
$$
\pr{S_n > na} \sim n\exp{\{-R(y_n)\}},\quad n\to\infty.
$$
Also, for any sequence $\varepsilon_n=o(\sqrt n)$,
$$
\pr{S_n>na}\sim\pr{S_n>na+\varepsilon_n}.
$$

\end{cor}

{\sc Proof.} This is just a reformulation of Lemma \ref{lem3}. We
should just note that in this case $R''(y_n)\sim 1/n$.

\begin{lem}\label{iterations}
Under the conditions of Theorem \ref{rw_greater_half} let
$t_n=na$. Define a sequence
$$y_n^{(0)}=na,\quad y_n^{(j)}=y_n^{(j-1)}-nR'(y_n^{(j-1)}).$$
Then, for any $j\ge 1/(2k)$,
$$
\pr{S_n>na}=n \exp\{-R(y_n^{(j)})\}.
$$
\end{lem}
{\sc Proof of lemma \ref{iterations}.} We have
\begin{eqnarray*}
|y_n^{(1)}-y_n^{(0)}|=n|R'(na)|=ng'(na).
\end{eqnarray*}
This implies that $y_n^{(2)}\sim na$. Assume that we proved
$y_n^{(i)}\sim na$ for all  $i<j$. Then, using regular variation
of $g''$ we obtain,
\begin{multline} \label{sub_1}
|y_n^{(j)}-y_n^{(j-1)}|=n|R'(y_n^{j-1})-R'(y_n^{j-2})|\\=(1+o(1))n|g'(y_n^{j-1})-g'(y_n^{j-2})|
=(1+o(1))n|g''(na)||y_n^{j-1}-y_n^{j-2}|=o(n).
\end{multline}
Therefore, we can argue by induction that for $j\ge 1$,
$$|y_n^{(j)}-y_n^{(j-1)}|=O(1)(n|g''(na)|)^{j-1}ng'(na)=
O(1)n(g'(n))^{j}=O(1)n\left(\frac{g(n)}{n}\right)^{j}.$$ Now make
use of Condition (\ref{def_k}), then
\begin{equation} \label{sub_2}
|y_n^{(j+1)}-y_n^{(j)}|=o(1)n
\left(\frac{n^{1-1/(k+2)}}{n}\right)^{j+1}=o(1)n^{1-(j+1)/(k+2)}=o(\sqrt
n),
\end{equation}
provided $j\ge k/2$. Then
$$R'(y_n^{(j)})=\frac{y_n^{(j)}-y_n^{(j+1)}}{n}=\frac{o(\sqrt n)}{n}=o(1/\sqrt n).$$
The proof of Lemma \ref{iterations} is complete.

\begin{lem} \label{lemma_subex}
Assume that all conditions of the previous lemma hold. Then the
sequence $a_n = \dfrac{\pr{S_n \ge 0}}{n}$ is subexponential.
\end{lem}

{\sc Proof of Lemma \ref{lemma_subex}.}

It follows from lemma \ref{iterations} that $a_n \sim
e^{-R(y_n^{(j)})}$ for sufficiently large $j$. We shall prove that
for any $j$,
\begin{equation} \label{sub_3}
\lim\limits_{n \to \infty} \frac{n R'(y_n^{(j)})}{R(y_n^{(j)})} =
\beta.
\end{equation}
This will imply that the sequence $a_n$ is subexponential due to
the sufficient conditions given in \cite{Kl}. We prove
(\ref{sub_3}) by induction. If $j=0$, then $R((y_n^{(j)}) =
g(y_n^{(j)})$ and (\ref{sub_3}) holds since $g$ is a regularly
varying function with parameter $\beta$. Assume (\ref{sub_3})
holds for some $j$. Note that (\ref{sub_2}) implies that
$R'(y_n^{(j)}) = \dfrac{y_n^{(j)} - y_n^{(j+1)}}{n}$ is regularly
varying, and hence, taking (\ref{sub_3}) into account,
$R(y_n^{(j)})$ is also regularly varying. Recall that $y_n^{(j)}
\sim na$ for each $j$ and that $y_n^{(j)} - y_n^{(j+1)} = o(n)$
(see (\ref{sub_1})). Then $R(y_n^{(j+1)}) = R(y_n^{(j)} +
(y_n^{(j+1)} - y_n^{(j)})) \sim R(y_n^{(j)})$ and also
$R'(y_n^{(j+1)}) \sim R'(y_n^{(j)})$. Hence, (\ref{sub_3}) holds
for $j+1$ as well. The proof of Lemma \ref{lemma_subex} is complete.

We now give the result on the asymptotic behaviour of the tails of
$\tau_x$ and $\nu_x$.

\begin{thm}\label{st_greater_half}
Assume that $\{X_t\}$ is a L\'evy process or a random walk such that
$\expect{X_1} = -a < 0$. Let all the conditions of Theorem
\ref{rw_greater_half} hold for the distribution of $X_1$. Let
$R'(y(t)) = o(1/\sqrt t)$. Then
\begin{eqnarray} \label{as_levy_greater_half1}
\pr{\tau_x>t}&\sim &\expect{\tau_x}\exp{\{-R(y(t))\}},\quad
t\to\infty,\\
\label{as_levy_greater_half2} \pr{\nu_x>n}& \sim
&\expect{\nu_x}\exp{\{-R(y(n))\}},\quad n\to\infty.
\end{eqnarray}
\end{thm}

{\sc Proof of Theorem \ref{st_greater_half}.} We use Theorem
\ref{thm_disc} again. First, it follows from Corollary~\ref{cor1}
that $\pr{X_n>y}\sim\pr{X_n>0}\sim\pr{X_{n+1}>0}$. Second,
according to lemma \ref{lemma_subex}, sequence
$\alpha_n=\pr{X_n>0}/n$ is subexponential. Then, we can just apply
Theorem \ref{thm_disc}. The proof of Theorem \ref{st_greater_half} is complete.

\begin{cor}\label{bp_greater_half}
Let $\pr{B_1>y}\sim e^{-g(y)}$. Assume that $g$ and $B_1$ satisfy
all conditions of Theorem \ref{st_greater_half} with $a=1-\rho$.
Then the asymptotics are given by
\begin{eqnarray} \label{as_bp_half1}
\pr{bp(x)>t}&\sim
&\frac{x}{\expect{A}(1-\rho)}\exp{\{-R(y(t))\}},\quad t\to\infty.
\end{eqnarray}
\end{cor}

{\sc Proof of corollary \ref{bp_greater_half}} repeats the proof
of corollary \ref{cor2}.

One can apply Lemma \ref{iterations} to obtain tail asymptotics
for the distributions of $\tau_x$ and $\nu_x$ explicitly. In
particular, the following corollary holds.

\begin{cor}
Under the conditions of Theorem \ref{st_greater_half} let
$g(y)=o(y^{3/4})$. Then
$$
\pr{\tau_x>t}\sim \expect{\tau_x}\exp{\{-R(ta-tg'(ta))\}},\quad
t\to\infty.
$$
If $\pr{B_1>t}\sim e^{-t^\beta},1/2<\beta<1$, then for
$k=[1/(1-\beta)]$ some positive constants $D_1,\ldots, D_k>0$,
$$
\pr{\tau_x>t}\sim\expect{\tau_x}\exp\{-(at)^\beta+D_1
t^{2\beta-1}+\ldots+D_kt^{k\beta-k+1}\}.
$$
\end{cor}
Similar corollaries can be formulated for $\nu_x$ and $bp(x)$.

\subsection{Cram\'er's case}

Let $m(s)=\expect{ e^{sX_1}}$ be the moment generating function of
$X_1$.
\begin{thm}\label{cramers_case}
Assume that $\{X_t\}$ is a L\'evy process or a random walk such that
$\expect{X_1} = -a < 0$. Let solution $\alpha$ to the equation
$m'(s)=0$ exist and $m''(\alpha)<\infty$. Put $\gamma=\ln
m(\alpha).$ Assume also that the distribution of $X_1$ is
non-lattice. Then
\begin{eqnarray*}
\pr{\tau_x>t}&\sim & V(x) \frac{1}{\sqrt{2\pi}t^{3/2}\widehat
\sigma(\alpha)\alpha}e^{-\gamma t},\quad t\to\infty,\\
\pr{\nu_x>n}&\sim &\frac{1}{\sqrt{2\pi}n^{3/2}\widehat
\sigma(\alpha)\alpha}e^{-\gamma n},\quad n\to\infty.
\end{eqnarray*}
\end{thm}

{\sc Proof of Theorem \ref{cramers_case}.} It follows from the
Petrov Theorem (see \cite[theorem 2]{petrov}) that $a_n=\pr{X_n\ge
0}/n\sim \frac{1}{\sqrt{2\pi}n^{3/2}\widehat
\sigma(\alpha)\alpha}e^{-\gamma n}.$ Then $a_n\in\mathcal
Ss(\gamma).$ Indeed, $a_{n-1}/a_n\to e^{\gamma}$ and for some
constant $C$,
$$
\sum_{k=1}^{n-1} \frac{a_k a_{n-k}}{a_n}\le C \sum_{k=1}^{n-1}
(\frac{n}{k(n-k)})^{3/2}\le 4C
\sum_{k=1}^{n-1}\frac{1}{k^{3/2}}<4C
\sum_{k=1}^{\infty}\frac{1}{k^{3/2}}.
$$  It follows from the dominated convergence Theorem that
$a_n^{*2}/a_n\to 2\sum_{n=1}^\infty e^{\gamma n}a_n$. Also, Petrov's
Theorem implies that
$$
\pr{X_n\ge y}/\pr{X_n\ge 0}\sim e^{-\alpha y},\quad n\to\infty.
$$
Therefore, we can apply Theorem \ref{thm_disc} to obtain the
needed asymptotics. The proof of Theorem \ref{cramers_case} is complete.

\begin{cor}\label{bp_cramers_case}
Let $\alpha>0$ be a solution to the equation $\lambda
m'_B(\alpha)=1$ such that $\widehat \sigma^2=\lambda
m''_B(\alpha)<\infty$. Put $\gamma=\alpha-\lambda(m_B(\alpha)-1)$.
Then,
\begin{equation}\label{asymp_bp_cramer}
\pr{bp(x)>t}\sim\frac{1}{\sqrt{2\pi\widehat \sigma^2}\gamma
t^{3/2}}xe^{\alpha x}e^{-\gamma t}
\end{equation}
\end{cor}

{\sc Proof of Corollary \ref{bp_cramers_case}.} It is clear that
$\alpha$ and $\gamma$ are exactly the same as in
Theorem~\ref{cramers_case}. Therefore, we should just find $V(x)$.
Since $e^{-\gamma}=\expect{e^{\alpha X_1}}$, the process
$\exp\{\alpha X_t+\gamma t\}$ is a martingale with mean $1$. Then,
since  $X_{\tau_x}=-x$, we have
$$
1=\expect{e^{\alpha X_{\tau_x}+\gamma \tau_x}}=e^{-\alpha
x}\expect{e^{\gamma \tau_x}}.
$$
Then,
$C(x)=\gamma^{-1}(\expect{e^{\gamma\tau_x}}-1)=\gamma^{-1}(e^{\alpha
x}-1)$, and it follows from (\ref{def_Vx}) that
$$
V(x)=C(x)+\alpha e^{\alpha x}\int_0^x e^{-\alpha y}C(y){\rm dy}=
\frac{\alpha}{\gamma}x e^{\alpha x}.
$$ The proof of Corollary \ref{bp_cramers_case} is complete.

\subsection{Intermediate case}

We now proceed to the intermediate case, that is when equation the
$m'(s)=0$ does not have a positive solution but $m(s)<\infty $ for
some $s>0$. In this case we shall assume that $\pr{\xi_1 > t} =
e^{-\alpha t} \overline G(t)$ for all positive $t$ where $\alpha >
0$ and $\overline G(t)$ is a tail of some heavy-tailed
distribution. Introduce the random walk $\{\tilde S_n\}$ (called
the adjunct random walk in \cite{bertoin_doney}) whose increments
have the distribution
\[
\tilde F(dy) = \frac{1}{m(\alpha)} e^{\alpha y} F(dy)
\]
and denote $\delta = - \expect{\tilde S_1}$.
The following result on large deviations may be found in
\cite{DDS}. It is a generalisation of the result of
\cite{bertoin_doney} where asymptotics for the large deviations
probabilities are obtained under the assumption that $\overline
G(t)$ is a regularly varying function.
\begin{thm} \label{ld_intermediate}
Assume that $\pr{\xi_1 > t} = e^{-\alpha t} \overline G(t)$ for
all positive $t$ where $\alpha > 0$ and $\overline G(t)$ satisfies
conditions of Theorem \ref{ld_less_half}. Let $\expect{\xi_1} <
0$, and $m'(s)\neq 0$ for $0<s\le \alpha$. Assume also that
$\delta < \infty$. Put $e^{-\gamma}=m(\alpha)$. Then uniformly in
$x$ such that $x \ge - n (\delta - \varepsilon)$
\[
\pr{S_n > x} \sim \frac{1}{m(\alpha)} e^{-\gamma n} e^{-\alpha x}
n \overline G(x + n \delta).
\]
\end{thm}

\begin{rem}
It is easy to see that the conditions of Theorem
\ref{ld_intermediate} imply that $\delta > 0$.
\end{rem}

Using the latter Theorem, one can obtain the following result for
the tail asymptotics of $\tau_x$ and $\nu_x$.

\begin{thm}\label{thm_intermediate}
Assume that $\{X_t\}$ is a L\'evy process or a random walk such
that the distribution of $X_1$ satisfies 
the conditions of Theorem \ref{ld_intermediate}. Then
\[
\pr{\tau_x>t} \sim  V(x) \frac{1}{m(\alpha)} e^{-\gamma t}
\overline G(t \delta)
\]
and
\[
\pr{\nu_x > n} \sim V_{rw}(x) \frac{1}{m(\alpha)} e^{-\gamma n}
\overline G(n \delta).
\]
\end{thm}

{\sc Proof of Theorem \ref{thm_intermediate}.}

Theorem \ref{ld_intermediate} implies that
\[
\pr{X_n \ge 0} \sim \frac{1}{m(\alpha)} e^{-\gamma n} n \overline
G(n \delta)
\]
and
\[
\pr{X_n \ge y} \sim \frac{1}{m(\alpha)} e^{-\gamma n} n \overline
G(n \delta + y) e^{-y \alpha}.
\]
The conditions of Theorem \ref{thm_disc} can now be checked in a
straightforward manner. The proof of Theorem \ref{thm_intermediate} is complete.

\begin{rem}
To the best of our knowledge, the only result on large deviations of
sums of random variables belonging to the intermediate case is
contained in \cite{bertoin_doney}. The result presented in this
paper concerns distributions $F$ such that the function $e^{\alpha
x}\overline F(x)$ is regularly varying with parameter $-\beta,
2<\beta<\infty$. The analogue of Theorem \ref{thm_intermediate} for
such distributions may be obtained using Lemma 3 from
\cite{bertoin_doney} instead of our Theorem \ref{ld_intermediate}.
\end{rem}

\appendix

\section{Proofs of Theorems \ref{thm1} and \ref{thm_disc}}

{\sc Proof of Remark \ref{equality}.} Indeed, prove first that $E
e^{\alpha X_1} \le e^{-\gamma}$. Fix arbitrary $C > 0$. Then for
large enough $t$, uniformly in $y \in (-C,C)$,
\[
e^{\alpha y} \le (1+\varepsilon) \frac{P(X_t > -y)}{P(X_t > 0)}.
\]
Consider now
\begin{multline*}
\int_{-C}^C e^{\alpha y} P(X_1 \in dy) \le (1+\varepsilon)
\frac{\int_{-C}^C P(X_t > -y) P(X_1 \in dy)}{P(X_t > 0)} \\ =
(1+\varepsilon) \frac{\int_{-C}^C P(X_t > -y) P(X_{t+1}-X_t \in
dy)}{P(X_t > 0)} \\ = (1+\varepsilon) \frac{P(X_{t+1}-X_t \in
(-C,C), X_{t+1}
> 0)}{P(X_t > 0)} \le \frac{P(X_{t+1} > 0)}{P(X_t > 0)} \le
(1+\varepsilon)^2 e^{-\gamma}
\end{multline*}
for large enough $t$. Since $\varepsilon$ is an arbitrary positive
number, we have $\int_{-C}^C e^{\alpha y} P(X_1 \in dy) \le
e^{-\gamma}$ and hence, $E e^{\alpha X_1} \le e^{-\gamma}$.

The inequality $E e^{\alpha X_1} \ge e^{-\gamma}$ can be proved in a
similar way. Take arbitrary $C>0$ and consider
$$
\int_{-C}^C \frac{P(X_t > -y) P(X_{t+1}-X_t \in dy)}{P(X_t > 0)} \le
(1+\varepsilon) \int_{-C}^C e^{\alpha y} P(X_1 \in dy) \le
(1+\varepsilon) E e^{\alpha X_1}
$$
for sufficiently large $t$. Since $C$ is arbitrary, we have
\[
\frac{P(X_{t+1} > 0)}{P(X_t > 0)} = \int_{-\infty}^\infty
\frac{P(X_t
> -y) P(X_{t+1}-X_t \in dy)}{P(X_t > 0)} \le (1+\varepsilon) E e^{\alpha
X_1}.
\]
We also have that for sufficiently large $t$
\[
e^{-\gamma} \le (1+\varepsilon)\frac{P(X_{t+1} > 0)}{P(X_t > 0)} \le
(1+\varepsilon)^2 E e^{\alpha X_1}
\]
which concludes the proof since the LHS does not depend on
$\varepsilon$. The proof of Remark \ref{equality} is complete.

{\sc Proof of Theorem \ref{thm1}.} Recall that $N_t=\inf_{s\le t}
X_s.$ It is clear that $\pr{\tau_x>t}=\pr{|N_t|\le x}$. Our starting
point is the formula which follows from the Wiener-Hopf identity for
L\'evy processes, see \cite[(47.9)]{Sato} (with obvious changes: we
should substitute the infimum process instead of the supremum
process). Thus, for $q>0$ and $u\ge 0$,
\begin{eqnarray}\label{wh}
q\int_0^\infty e^{-qt} \expect{e^{uN_t}}{\rm dt}=
\exp\left\{\int_0^\infty t^{-1}e^{-qt} {\rm dt}\int_{(-\infty, 0)}
(e^{uy}-1) \pr{X_t\in\rm{dy}}\right\}.
\end{eqnarray}
For $q>0$, make use of the Frullani integral
$$
-\ln q=\int_0^\infty \frac{e^{-qt}-e^{-t}}{t}\rm dt
$$
and rewrite (\ref{wh}) in the following form
\begin{multline}\label{wh2}
\int_0^\infty e^{-qt} \expect{e^{uN_t}}{\rm dt}=
\exp\biggl\{\int_0^\infty \frac{e^{-qt}-e^{-t}}{t}{\rm dt}\\
+\int_0^\infty t^{-1}e^{-qt} {\rm dt}\int_{(-\infty, 0)}
(e^{uy}-1)\pr{X_t\in\rm{dy}}\biggr\}.
\end{multline}
First we want to show that the right-hand side in (\ref{wh2})
converges as $q\to 0$. For that, let us represent the exponent on
the right-hand side of (\ref{wh2}) as follows,
\begin{multline}\label{eq1}
\int_0^1 \frac{e^{-qt}-1}{t}{\rm dt}+\int_0^1
\frac{1-e^{-t}}{t}{\rm dt}-\int_1^{\infty}\frac{e^{-t}}{t}{\rm
dt}\\+
\int_1^{\infty}t^{-1}e^{-qt}\left\{1+\int_{(-\infty, 0)} (e^{uy}-1)\mathbf P(X_t\in\rm{dy})\right\}{\rm dt}\\
+\int_0^1 t^{-1}e^{-qt} {\rm dt}\int_{(-\infty, 0)}
(e^{uy}-1)\mathbf P(X_t\in\rm{dy}).
\end{multline}
It is clear that the first integral converges to 0 as $q\to 0$.
The second and third integral are well defined and constant. The
fourth and fifth integral are monotone in $q$, and therefore it is
sufficient to prove that they are finite for $q=0$ and then to
apply the monotone convergence Theorem. For $q=0$, the fourth
integral is equal to
\begin{multline}\label{fourth_integral}
\int_1^{\infty}t^{-1}\left\{1+\int_{(-\infty, 0)} (e^{uy}-1)\pr{X_t\in\rm{dy}}\right\}{\rm dt}\\
=\int_1^{\infty}t^{-1}\left\{\pr {X_t\ge 0}+\int_{(-\infty, 0)}
e^{uy}\pr{X_t\in\rm{dy}}\right\}{\rm dt}.
\end{multline}
To deal with the second term in (\ref{fourth_integral}), we prove
\begin{lem}\label{lem1}
Let function $\pr{X_t\ge 0}/t$ be such that Condition (\ref{lt})
holds.
  Assume also that condition (\ref{cond1}) of Theorem \ref{thm1}
holds.
Then, there exists $u_0$ such that for any $u>u_0,$
\begin{eqnarray}\label{eq12}
\int_{(-\infty, 0)} e^{uy}\pr{X_t\in{\rm dy}}\sim
\frac{\alpha}{u-\alpha} \pr{X_t\ge 0}.
\end{eqnarray}

\end{lem}
{\bf Remark.} We use the convention that $a(x)\sim 0\cdot b(x)$
means that $a(x)=o(b(x)).$

 {\sc Proof of Lemma \ref{lem1}.}
Using integration by parts we obtain that
$$
\int_{(-\infty, 0)} e^{uy}\pr{X_t\in{\rm dy}}=u\int_0^\infty
(\pr{X_t\ge-y}-\pr{X_t\ge 0})e^{-uy}{\rm dy}.
$$
We can pick function $h(t)\uparrow\infty$ such that
$(\ref{cond1})$ holds uniformly in $y\in[-h(t), 0]$. Then, for
$u>\alpha$,
\begin{multline*}
u \int_0^{h(t)}
(\pr{X_t\ge-y}-\pr{X_t\ge 0})e^{-uy}{\rm dy}\\
=(1+o(1))\pr{X_t\ge 0} u \int_0^{h(t)} (e^{\alpha y}-1)e^{-uy}{\rm
dy}= \frac{\alpha+o(1)}{u-\alpha}\pr{X_t\ge 0}
\end{multline*}

Further, since $\pr{X_t\ge 0}/t$ satisfies Condition (\ref{lt}),
for sufficiently large $t$, \begin{eqnarray}\label{eq11}
\pr{X_{t+n}\ge 0}\le e^{n}\pr{X_t\ge 0}, \quad n\ge 1.
\end{eqnarray}
Since $\pr{X_1> 0}>0$, then for some $\delta=1/l>0$,
 $\pr{X_1\ge \delta}>0$, where $l$ is a positive integer.
Then,
$$
\pr{X_{ln}\ge n}\ge\pr{X_1\ge \delta,X_2-X_1\ge \delta,\ldots
X_{ln}-X_{ln-1}\ge \delta}=\pr{X_1\ge\delta}^{ln}.
$$
Let $u_0$ be such that $e^{-u_0}=e^{-2l}\pr{X_1\ge\delta}^{l}$.
 Then,  for all $u>u_0$,
$$
e^{-un}\le e^{-u_0n}=e^{-2ln}\pr{X_1\ge \delta}^{ln}\le
e^{-2ln}\pr{X_{ln}\ge n}.
$$
Therefore, for $u>u_0$,
\begin{multline*}
\int_{-\infty}^{-h(t)}(\pr{X_t\ge-y}-\pr{X_t\ge 0})e^{-uy}{\rm dy}
\le  \sum_{n=[h(t)]}^\infty e^{-un}\pr{X_t\ge -n}\\
\le \sum_{n=[h(t)]}^\infty e^{-2ln}\pr{X_{ln}\ge n}\pr{X_t\ge -n}\\=
 \sum_{n=[h(t)]}^\infty e^{-2ln}\pr{X_{t+ln}-X_t\ge n,X_t\ge -n}
\le  \sum_{n=[h(t)]}^\infty e^{-2ln}\pr{X_{t+ln}\ge0}.
\end{multline*}
It follows from (\ref{eq11}) that, as $t\to\infty$,
$$
\sum_{n=[h(t)]}^\infty e^{-2ln}\pr{X_{t+ln}\ge 0}\le \pr{X_t\ge 0}
\sum_{n=[h(t)]}^\infty e^{-ln}=o(\pr{X_t\ge 0}).
$$

The proof of  Lemma \ref{lem1} is complete.

We may now continue to analyse (\ref{fourth_integral}). It follows
from Lemma \ref{lem1} that for some constant $C>0$,
$$
\int_1^{\infty}t^{-1}\left\{\pr {X_t\ge 0}+\int_{(-\infty, 0)}
e^{uy}\pr{X_t\in\rm{dy}}\right\}{\rm dt}\le C
\int_1^{\infty}t^{-1}\pr {X_t\ge 0}{\rm dt} <\infty.
$$
The finiteness of  the latter integral follows from (\ref{trans}).
 We now proceed to  the last term in (\ref{eq1}). Making use of
the inequality: $1-e^{-x}\le x, $ for all $x\ge 0$, we obtain
\begin{multline*}
\int_0^1 t^{-1}{\rm dt}\int_{(-\infty, 0)} (1-e^{uy})\mathbf
P(X_t\in\rm{dy})\\\le \int_0^1 t^{-1}{\rm dt} \left(\int_{(-1, 0)}
(-uy)\mathbf P(X_t\in{\rm{dy}})+\pr{X_t<-1}\right)\\
 =\int_0^1\expect{u X_t {\bf I}(X_t\in[-1,0]) }t^{-1}{\rm dt}+\int_0^1 t^{-1}{\rm dt} \pr{X_t<-1}<\infty.
\end{multline*}
The finiteness of the latter integral follows from the estimates
in Lemma \ref{lem_sato} in Appendix A.  Therefore,  the last term
in (\ref{eq1}) converges by the monotone convergence Theorem.

Now, letting $q\to 0$ in (\ref{wh2}), we have,
$$
\int_0^\infty \expect{e^{uN_t}}{\rm dt}<\infty.
$$

For fixed $u$, let now
$$
f_u(t)=\frac{\expect{e^{uN_t}}}{\int_0^\infty
\expect{e^{uN_t}}{\rm dt}}
$$
be the density of a random variable $Z$. Then, using
representations (\ref{eq1}) and (\ref{fourth_integral}), we may
rewrite (\ref{wh2})
\begin{multline*}
\expect{e^{-qZ}}=\exp\biggl\{\int_0^1 \frac{e^{-qt}-1}{t}{\rm dt}+
\int_1^{\infty}\frac{e^{-qt}-1}{t}\biggl( \pr
{X_t\ge 0}+\int_{(-\infty, 0)} e^{uy}\pr{X_t\in\rm{dy}}\biggr){\rm dt}\\
+\int_0^1 \frac{e^{-qt}-1}{t} {\rm dt}\int_{(-\infty, 0)}
(e^{uy}-1)\mathbf P(X_t\in\rm{dy})\biggr\}\\
\equiv \exp\biggl\{\int_0^1 ({e^{-qt}-1})\nu_1({\rm dt})+
\int_1^{\infty}({e^{-qt}-1})\nu_2({\rm dt}) +\int_0^1
({e^{-qt}-1})\nu_3({\rm dt})\biggr\}.
\end{multline*}
Then  $Z$ is an infinitely divisible variable on $[0,\infty)$ with
the L\'evy measure $\nu(dt)=\nu_1(dt)+\nu_2(dt)+\nu_3(dt)$. Indeed,
$\int_0^1 t \nu_1(dt)<\infty$. Further, as we have already shown
$\nu_3(0,1)<\infty$. Finally, as follows from Lemma~\ref{lem1},
$$
f_2(t)\equiv\frac{d\nu_2}{dt}=\frac{\pr {X_t\ge
0}}{t}+\frac{1}{t}\int_{(-\infty, 0)}
e^{uy}\pr{X_t\in{\rm{dy}}}\sim\frac{u}{u-\alpha}\frac{\pr {X_t\ge
0}}{t},\quad t\to\infty.
$$
Therefore, by Condition  (\ref{trans}), $\nu_2(1,\infty)<\infty.$
Now we are in the position to apply Theorem \ref{subex_density}
from the Appendix. Since the density of the L\'evy measure
$f(t)\sim\pr{X_t\ge 0}/t$ belongs to the class $\mathcal
Sd(\gamma)$, by Theorem \ref{subex_density}, we have for any fixed
$u$,
$$
f_u(t) \sim \left(\int_0^\infty e^{\gamma y}f_u(y){\rm
dy}\right)\frac{u}{u-\alpha}\frac{\pr{X_t\ge 0}}{t}.
$$
Equivalently, for all $u$,
\begin{eqnarray}\label{eq2}
\frac{\expect{e^{uN_t}}}{\pr{X_t\ge 0}/t}\to
\frac{u}{u-\alpha}\int_0^\infty e^{\gamma t}\expect{e^{uN_t}}{\rm
dt}.
\end{eqnarray}
Then, changing the order of integration, we obtain
$$
\int_0^\infty e^{\gamma t}\expect{e^{uN_t}}{\rm dt}= \int_0^\infty
e^{\gamma t}\left(\int_0^\infty {e^{-ux}}\pr{|N_t|\in {\rm
dx}}\right){\rm dt}=
\int_0^\infty e^{-ux}d_xC(x),
$$
where $C(x)=\int_0^\infty e^{\gamma t}\pr{|N_t|\le x} {\rm dt}$.
Therefore, (\ref{eq2}) is equivalent to
\begin{eqnarray*}
\frac{\expect{e^{uN_t}}}{\pr{X_t\ge 0}/t}\to
\frac{u}{u-\alpha}\int_0^\infty e^{-ux}d_xC(x).
\end{eqnarray*}
Now note that $u/(u-\alpha)$ is the LST  of the measure $D$ which
has unit mass at $0$ and density $\alpha e^{\alpha y}$ on the
positive half-line. Therefore,
\begin{eqnarray}
\frac{\pr{|N_t|\le x}}{\pr{X_t\ge 0}/t}\to D*C(x)&=&
V(x)=C(x)+\alpha e^{\alpha x} \int_0^x e^{-\alpha y}C(y){\rm
dy}\label{def_Vx}\\\nonumber &=&e^{\alpha x}\int_0^\infty
e^{\gamma t}\expect{e^{\alpha N_t};|N_t|\le x}{\rm dt}
\end{eqnarray}
for all $x$, where function $V(x)$ is continuous. This is
equivalent to
$$
\frac{\pr{\tau_x>t}}{\pr{X_t\ge 0}/t}\to V(x).
$$
Finally, it is clear that if $\alpha=\gamma=0$ then
$$
V(x)=\int_0^\infty \pr{|N_t|\le x}{\rm dt}= \int_0^\infty
\pr{\tau_x>t}{\rm dt}=\expect{\tau_x}.
$$
The proof of Theorem \ref{thm1} is complete.

Now we will show that for Theorem \ref{thm1} to hold it is
sufficient that its conditions hold for positive integers $t$.
This will allow us to reduce the problem of verifying properties
of $X_t,
 t\in{\mathbf R}$ to verifying the corresponding properties of the random walk $X_n, n\in\mathbf
N$.

\begin{lem} \label{lemma_asymptotics_levy}
Assume that sequence $\pr{X_n\ge 0}/n,n\in\mathbf N$ belongs
to the class $Ss(\gamma)$, 
for any fixed $y$
\begin{equation} \label{ltxt2}
\mathsf P(X_{n} \ge  0) \sim e^{\alpha y}\mathsf P(X_{n} \ge
y),\quad n\to \infty,
\end{equation}
and $e^{-\gamma}=\expect{e^{\alpha X_1}}.$ Then, function
$\pr{X_t\ge 0}/t, t\in\mathbf R$ belongs to the class
$\mathcal Sd(\gamma)$ 
and for any fixed $y$,
\begin{equation}\label{ltxt2_cont}
\mathsf P(X_{t} \ge  0) \sim e^{\alpha y}\mathsf P(X_{t} \ge
y),\quad t\in\mathbf R^+,t\to\infty.
\end{equation}
\end{lem}

{\it Proof of Lemma \ref{lemma_asymptotics_levy}.} First, we prove
that for any $0<\varepsilon< 1$, \begin{eqnarray}\label{eq15}
\pr{X_{n+\varepsilon} \ge 0}\sim e^{-\gamma\varepsilon}\pr{X_n \ge
0}.
\end{eqnarray}
It is not difficult to prove that there exists a function
$h(n)\uparrow\infty$ such that Condition (\ref{ltxt2}) holds
uniformly in $|z|\le h(n).$ We start with the total probability
formula
\begin{multline}\label{tot_probab}
\pr{X_{n+\varepsilon} \ge 0}\equiv P_1+P_2+P_3=\pr{X_{n+\varepsilon} \ge 0,|X_{n+\varepsilon}-X_n|\le h(n)}\\
+\pr{X_{n+\varepsilon} \ge
0,X_{n+\varepsilon}-X_n>h(n)}+\pr{X_{n+\varepsilon} \ge
0,X_{n+\varepsilon}-X_n<-h(n)}.
\end{multline}
Then, since $\expect {e^{\alpha X_1}}=e^{-\gamma}$,
$$
P_1=\int_{-h(n)}^{h(n)}\pr{X_\varepsilon\in{\rm dy}}\pr{X_n\ge
-y}\sim \pr{X_n\ge 0} \int_{-\infty}^\infty e^{\alpha
y}\pr{X_\varepsilon\in{\rm dy}}=e^{-\gamma\varepsilon}\pr{X_n\ge
0}.
$$
Before proceeding further, note that if we take $\varepsilon=1$
then it follows from (\ref{lt_disc}), (\ref{tot_probab}) and the
latter equivalence that
\begin{eqnarray}\label{eq16}
\pr{X_{n+1} \ge 0,X_{n+1}-X_n>h(n)}=o(\pr{X_n\ge 0}).
\end{eqnarray}
Further,
$$
P_2=\pr{X_n\ge-h(n)}\pr{X_{n+\varepsilon}-X_n>h(n)}+\int_{-\infty}^{-h(n)}
\pr{X_n\in {\rm dy}}\pr{X_{n+\varepsilon}-X_n\ge-y}.
$$
Now note that
$$
\pr{X_1\ge y}\ge \pr{X_\varepsilon\ge y,X_1-X_\varepsilon\ge 0}=
\pr{X_\varepsilon\ge y}\pr{X_{1-\varepsilon}\ge 0},
$$
which implies that
\begin{multline*}
P_2\le \frac{1}{\pr{X_{1-\varepsilon}\ge 0}}(\pr{X_n\ge
-h(n)}\pr{X_{n+1}-X_n>h(n)}\\+\int_{-\infty}^{-h(n)} \pr{X_n\in
{\rm dy}}\pr{X_{n+1}-X_n\ge -y}) =\frac{\pr{X_{n+1} \ge
0,X_{n+1}-X_n\ge h(n)}}{\pr{X_{1-\varepsilon}\ge 0}}.
\end{multline*}
After applying (\ref{eq16}) it is clear that $P_2=o(\pr{X_n>0})$.
Finally,
\begin{multline*}
P_3\le \pr{X_n\ge h(n),X_{n+\varepsilon}-X_n<-h(n)}\\\le
\pr{X_n\ge 0}\pr{X_{\varepsilon}<-h(n)}=o(\pr{X_n\ge 0}).
\end{multline*}

Now we should make use of the fact that if (\ref{eq15}) holds for
any fixed $\varepsilon\in(0,1)$, then, it holds uniformly in
$\varepsilon\in(0,1)$. Consequently, (\ref{lt}) holds for
$a(t)=\pr{X_t\ge 0}/t$. The proof of (\ref{ltxt2_cont}) is
similar. Finally, condition (\ref{subex}) for function $a(t)$
follows from the dominated convergence Theorem and the fact that
for some constant $C$,
$$
\int_1^{y-1}\frac{a(y-t)a(t)}{a(t)}{\rm dt }\le
C\sum_{k=1}^{[y-1]}\frac{a(k)a([y-1]-k)}{a[y-1]}<\infty.
$$
The proof of Lemma \ref{lemma_asymptotics_levy} is complete.

{\sc Proof of Theorem \ref{thm_disc}.} For L\'evy processes the
result follows directly from Theorem \ref{thm1} and Lemma
\ref{lemma_asymptotics_levy}. For random walks it can be proved
along the lines of \cite{doney}. The only difference is that we
should apply our Lemma \ref{lem1} instead of Lemma 4 in
\cite{doney}.

\section{L\'evy processes}

In this section we collect some facts from the theory of L\'evy
processes that we use in this paper.

This Lemma and its proof may be found in \cite[Lemma 30.3]{Sato}.
\begin{lem}\label{lem_sato}
Let $X_t$ be a L\'evy process on $\bf R^d$. For any $\varepsilon>0$
there is $C=C(\varepsilon)$ such that, for any $t$,
\begin{eqnarray*}
 \pr{|X_t|>\varepsilon}\le Ct.
\end{eqnarray*}
There are $C_1, C_2$ and $C_3$ such that, for any $t$ ,
\begin{eqnarray*}
\expect{|X_t|^2;|X_t|\le 1}&\le & C_1 t,\\
 |\expect{X_t;|X_t|\le
1}|&\le & C_2 t,\\
\expect{|X_t|;|X_t|\le 1}&\le& C_3 t^{1/2}.
\end{eqnarray*}

\end{lem}

The next Theorem is a version of \cite[theorem 25.3]{Sato} adapted
to our needs.
\begin{thm}\label{exp_bound}
Let $X_t$ be a L\'evy process on $\bf R^d$ with the L\'evy measure
$\nu$. Then, $X_t$ has a finite exponential moment if and only if
$[\nu]_{|x|>1}$ has a finite exponential moment.
\end{thm}

Let $F$ be an infinitely divisible law on $[0,+\infty).$ Its
Laplace law can be expressed as (\cite[theorem 30.1]{Sato})
$$
\int_0^\infty e^{-\lambda x} F(dx)=\exp\left\{-\gamma
\lambda-\int_0^\infty (1-e^{-\lambda x})\nu({\rm dx}))\right\},
$$
where $\gamma\ge 0$ is a constant and $\nu$ is a Borel measure on
$(0,\infty)$ for which $\mu\equiv\nu(1,\infty)<\infty$ and
$\int_0^1x\nu({\rm dx})<\infty.$

\begin{thm}\label{thm_egv} {\rm (\cite[theorem 1]{egv})}
For $F$ infinitely divisible on $[0,+\infty)$, the following
assertions are equivalent:

{\rm (i)}   $F\in\mathcal S$;

{\rm (ii)}  $\mu^{-1}\nu(1,x]\in\mathcal S$;

{\rm (iii)} $\overline F(x)\sim\overline \nu(x)$.

\end{thm}

We also need a density version of this Theorem.

\begin{thm}\label{subex_density}
Let the infinitely divisible law $F$ has a density $f$. Assume
that there is $x_0$ such that $\nu$ has a density $g(x)$ for
$x>x_0$. If $g(x)$ belongs to the class $\mathcal Sd(\gamma)$,
then

$$\lim_{x\to\infty}\frac{f(x)}{g(x)}=\int_0^\infty e^{\gamma x}f(x){\rm dx}.$$

\end{thm}
One can prove this Theorem exactly like \cite[theorem 3.1]{pakes}
for a distribution function from $\mathcal S(\gamma)$.
Corresponding properties of the class $\mathcal Sd(\gamma)$ may be
found in \cite[Section 3]{Kl2}.

{\bf Acknowledgments.}
The research of Denis Denisov was supported by the Dutch BSIK project ({\it BRICKS}) and
the EURO-NGI project. The research of V. Shneer was  supported by a James-Watt Scholarship from the Heriot-Watt
University, a United Kingdom Overseas Research Scholarship and the
EURO-NGI project.

We would like to thank Serguei Foss for
drawing attention to this problem and a number of helpful
discussions. We are also grateful to Onno Boxma, Ton Dieker and Bert
Zwart for many useful comments and suggestions that helped to
improve the manuscript. A part of this work was done while Vsevolod
Shneer was visiting EURANDOM.  He would like to thank EURANDOM for
the hospitality.

\end{document}